\documentclass[12pt]{article}
\usepackage{amssymb,amsmath,amsthm}
\usepackage{bm}
\usepackage{indentfirst}
\usepackage{amsfonts}
\usepackage{authblk}
\usepackage{times}
\usepackage{sectsty}
\usepackage{titlesec}
\catcode240=13 \def ð{\u{g}} \catcode231=13 \def ç{\c{c}} \catcode246=13 \def ö{\"{o}} \catcode254=13 \def þ{\c{s}}
\catcode252=13 \def ü{\"{u}} \catcode253=13 \def ý{{\i}} \catcode221=13 \def Ý{\.{I}} \catcode199=13 \def Ç{\c{C}}
\catcode208=13 \def Ð{\u{G}} \catcode214=13 \def Ö{\"{O}} \catcode222=13 \def Þ{\c{S}} \catcode220=13 \def Ü{\"{U}}
\providecommand{\keywords}[1]{\textbf{\textit{Keywords:}} #1}
\usepackage[a4paper,left=3cm,right=3cm,top=3cm,bottom=3cm]{geometry}
\usepackage{setspace}
\theoremstyle{definition}

\numberwithin{equation}{section}

\newtheorem{lem}{Lemma}

\newtheorem{teo}{Theorem}

\newtheorem{cor}{Corollary}

\newtheorem{ex}{Example}

\newtheorem*{pr}{Proof}
\newtheorem*{dem}{Demonstration}

\begin{document}

\title{Identities for the $\mathcal{L}_n$-transform, the $\mathcal{L}_{2n}$-transform and the $\mathcal{P}_{2n}$ transform and their applications}

\author[1  \footnote{Correspondence: ndernek@marmara.edu.tr}]{Nese Dernek}
\author[2]{Fatih Aylikci}

%\author[2]{Corresponding Author\thanks{email@2nduniversity.com}}
\affil[1]{Department of Mathematics, Marmara University, Istanbul, Turkey}
\affil[2]{Department of Mathematical Engineering, Yildiz Technical University, Istanbul, Turkey}
%\author{Neþe Dernek, \\ Department of Mathematics, University of Marmara, Istanbul, Turkey \and Fatih Aylýkcý \\ Department of Mathematics}

\maketitle

\begin{abstract}
In the present paper, the authors introduce several new integral transforms inclu-\\ding the $\mathcal{L}_n$-transform, the $\mathcal{L}_{2n}$-transform and $\mathcal{P}_{2n}$-transform generalizations of the classical Laplace transform and the classical Stieltjes transform as respectively. It is shown that the second iterate of the $\mathcal{L}_{2n}$-transform is essentially the $\mathcal{P}_{2n}$-transform. Using this relationship, a few new Parseval-Goldstein type identities are obtained. The theorem and the lemmas that are proven in this article are new useful relations for evaluating infinite integrals of special functions. Some related illustrative examples are also given.
\end{abstract}


%\begin{keywords}
%$L_2$-transform, Laplace transform, $L_n$-transform, Widder potential transform, Stieltjes transform, Parseval-Goldstein type theorem
%\end{keywords}

\section{Introduction, definitions}

Integral transforms are used in several branches of applied mathematics. Among these the Laplace transform is the most commonly used one in applications. \\
\indent Widder \cite{widder} presented a systematic account of the so-called Widder potential transform:
\begin{equation}
\mathcal{P}\{f(x);y\}=\int\limits_0^\infty \frac{xf(x)}{x^2 + y^2}dx. \label{widder}
\end{equation}

Widder pointed out that the potential transform is related to the Poisson integral re-\\presentation of a function which is harmonic in a half plane and give several inversion formula for the transform.\\
\indent The classic Laplace transform is defined as
\begin{equation}
\mathcal{L}\{f(x);y\}=\int\limits_0^\infty \exp(-xy) f(x) dx.
\end{equation}
Goldstein \cite{goldstein} introduced the following Parseval-type theorem for the classical Laplace transform:
\begin{equation}
\int\limits_0^\infty f(x) \mathcal{L}\{g(y);x\}dx=\int\limits_0^\infty g(y) \mathcal{L}\{f(x);y\}dy.
\end{equation}
Yurekli \cite{osman2} gave the following Parseval-Goldstein type theorem,
\begin{equation}
\int\limits_0^\infty \mathcal{L}\{f(x);y\}\mathcal{L}\{g(y);x\}dx=\int\limits_0^\infty g(y) \mathcal{S}\{f(x);y\} dy,
\end{equation}
involving the Laplace transform and the Stieltjes transform,
\begin{equation}
\mathcal{S}\{f(x);y\}=\int\limits_0^\infty \frac{f(x)}{x+y} dx. \label{stieltjes}
\end{equation}
Srivastava and Singh \cite{srivastava} gave the following Parseval-Goldstein type formula for the Widder potential transform:
\begin{equation}
\int\limits_0^\infty y \mathcal{P}\{f(x);y\} g(y) dy = \int\limits_0^\infty x \mathcal{P}\{g(y);x\} f(x) dx.
\end{equation}
\indent Yurekli and Sadek \cite{osman1} presented a systematic account of so called the $\mathcal{L}_2$-transform,
\begin{equation}
\mathcal{L}_2\{f(x);y\}= \int\limits_0^\infty x \exp(-x^2 y^2) f(x) dx.
\end{equation}
The $\mathcal{L}_2$-transform is related to the classical Laplace transform, via the identities
\begin{equation}
\mathcal{L}_2\{f(x);y\} = \frac{1}{2} \mathcal{L}\{f(\sqrt{x});y^2\},
\end{equation}
\begin{equation}
\mathcal{L}\{f(x);y\}=2\mathcal{L}_2\{f(x^2);\sqrt{y}\}.
\end{equation}
\indent Dernek, Srivastava and Yurekli \cite{dernek1} presented the $\mathcal{L}_4$-transform and the $\mathcal{P}_4$-transform respectively:
\begin{equation}
\mathcal{L}_4\{f(x);y\}=\int\limits_0^\infty x^3 \exp(-x^4 y^4)f(x)dx,
\end{equation}
\begin{equation}
\mathcal{P}_4\{f(x);y\}=\int\limits_0^\infty \frac{x^3 f(x)}{x^4 + y^4} dx.
\end{equation}
The $\mathcal{L}_4$-transform is related to the Laplace transform and the $\mathcal{L}_2$-transform by means of the following identities:
\begin{equation}
\mathcal{L}_4\{f(x);y\}=\frac{1}{4} \mathcal{L}\{f(x^{\frac{1}{4}});y^4\},
\end{equation}
\begin{equation}
\mathcal{L}_4\{f(x);y\}=\frac{1}{2}\mathcal{L}_2\{f(x^{\frac{1}{2}});y^2\}.
\end{equation}

\noindent Dernek, Srivastava and Yurekli \cite{dernek1,dernek2} gave the Parseval-Goldstein type relation between the $\mathcal{L}_4$-transform and the $\mathcal{P}_4$-transform as follows:
\begin{equation}
\int\limits_0^\infty x^3 f(x) \mathcal{P}_4\{g(u);x\} dx = \int\limits_0^\infty u^3 g(u) \mathcal{P}_4\{f(x);u\} du.
\end{equation}
\indent We now introduce a few generalized integral transforms.
The $\mathcal{L}_n ~(n=2^k, k=0,1,2,\ldots, \\ k \in \mathbb{N})$ transform defined in by the following formula,
%\indent We now introduce a few generalized integral transforms. The $\mathcal{L}_n$ $(n=2^k, k=0,1,2,...)$ transform and the $\mathcal{L}_{2n}$ transform as a generalization of the Laplace transform and the $\mathcal{L}_2$-transform are defined as
\begin{equation}
\mathcal{L}_n\{f(x);y\}=\int\limits_0^\infty x^{n-1}\exp(-y^n x^n) f(x) dx  ,~~~~ n=2^k,~~k=0,1,2,..., \label{laplacen}
\end{equation}
which is a generalization of the Laplace transform. \\
 The $\mathcal{L}_{2n}$ transform as a generalization of the $\mathcal{L}_2$-transform is defined as
\begin{equation}
\mathcal{L}_{2n}\{f(x);y\}=\int\limits_0^\infty y^{2n-1}\exp(-y^{2n} x^{2n}) f(x) dx, ~~~~ n \in \mathbb{N}. \label{laplace2n}
\end{equation}
The $\mathcal{L}_n$-transform and the $\mathcal{L}_{2n}$-transform are related to the Laplace transform and the $\mathcal{L}_2$-transform with the following formulas:
\begin{equation}
\mathcal{L}_n\{f(x);y\}=\frac{1}{n} \mathcal{L}\{f(x^{\frac{1}{n}});y^n\} ~, n=2^k, k \in \mathbb{N},  \label{laplacen.1}
\end{equation}
\begin{equation}
\mathcal{L}_{2n}\{f(x);y\}=\frac{1}{2n} \mathcal{L}\{f(x^{\frac{1}{2n}});y^{2n}\}, ~n \in \mathbb{N} \label{laplace2n.2}
\end{equation}
and
\begin{equation}
\mathcal{L}_n\{f(x);y\}= \frac{2}{n}\mathcal{L}_2\{f(x^{\frac{2}{n}});y^{\frac{n}{2}}\},
\end{equation}
\begin{equation}
\mathcal{L}_{2n}\{f(x);y\}=\frac{1}{n} \mathcal{L}_2\{f(x^{\frac{1}{n}});y^n\}. \label{laplace2n.1}
\end{equation}
\indent We introduce the $\mathcal{P}_n$-transform and the $\mathcal{P}_{2n}$-transform as follows:
\begin{equation}
\mathcal{P}_n\{f(x);y\}=\int\limits_0^\infty \frac{x^{n-1}f(x)}{x^n + y^n}dx,~~~~ n=2^k,~k=0,1,2,... ;
\end{equation}
\begin{equation}
\mathcal{P}_{2n}\{f(x);y\}=\int\limits_0^\infty \frac{x^{2n-1}f(x)}{x^{2n} + y^{2n}}dx,~~~~ n \in \mathbb{N}. \label{potential}
\end{equation}
The $\mathcal{P}_n$-transform is related the Stieltjes-transform  (\ref{stieltjes}) by means of the identity
\begin{equation}
\mathcal{P}_n\{f(x);y\}=\frac{1}{n} \mathcal{S}\{f(x^{\frac{1}{n}});y^n\}
\end{equation}
and the $\mathcal{P}_n$-transform is related the Widder potential transform (\ref{widder}) by means of the identity
\begin{equation}
\mathcal{P}_n\{f(x);y\}=\frac{2}{n} \mathcal{P}\{f(x^{\frac{2}{n}});y^{\frac{n}{2}}\}.
\end{equation}
Similarly the $\mathcal{P}_{2n}$-transform is related the Widder potential transform (\ref{widder}) by means of the identity
\begin{equation}
\mathcal{P}_{2n}\{f(x);y\}=\frac{1}{n} \mathcal{P}\{f(x^{\frac{1}{n}});y^n\}.
\end{equation}

\section{Theorems and corollaries}
\begin{lem}
The following identity holds true
\begin{equation}
\mathcal{L}_n\{\mathcal{L}_{2n}\{f(x);u\};y\}=\frac{\sqrt{\pi}}{2n} \int\limits_0^\infty x^{n-1} f(x) \exp(\frac{1}{4}\frac{y^{2n}}{x^{2n}})erfc(\frac{1}{2}\frac{y^n}{x^n})dx,~~~~n=2^k,~k=0,1,2,... \label{lemma1}
\end{equation}
provided that the integrals involved converge absolutely.
\pr
By the definitions (\ref{laplacen}) and (\ref{laplace2n}) we have
\begin{equation}
\mathcal{L}_n\{\mathcal{L}_{2n}\{f(x);u\};y\}=\int\limits_0^\infty u^{n-1} \exp(-u^n y^n)\Big[\int\limits_0^\infty x^{2n-1} \exp(-u^{2n} x^{2n}) f(x) dx\Big] du.
\end{equation}
With changing the order of integration which is permissible by absolute convergence of the integrals involved and using the definition of ~~$erfc(x)$~~ function we find
\begin{equation}
\mathcal{L}_n\{\mathcal{L}_{2n}\{f(x);u\};y\}=\int\limits_0^\infty x^{2n-1} f(x) \Big[\int\limits_0^\infty u^{n-1} \exp(-u^{2n} x^{2n} - y^n u^n )du \Big] dx \nonumber
\end{equation}
\begin{equation}
=\int\limits_0^\infty x^{2n-1} f(x) \exp(\frac{1}{4}\frac{y^{2n}}{x^{2n}})\Big[\int\limits_0^\infty u^{n-1} \exp(-x^{2n}(u^n + \frac{1}{2}\frac{y^n}{x^{2n}})^2) du\Big]dx.
\end{equation}
Making the following change of variable,
\begin{equation}
x^n(u^n + \frac{1}{2}\frac{y^n}{x^{2n}})=t,
\end{equation}
we obtain the assertion (\ref{lemma1}).
%\begin{equation}
%\mathcal{L}_n\{\mathcal{L}_{2n}\{f(x);u\};y\}=\frac{\sqrt{\pi}}{2n} \int\limits_0^\infty x^{n-1} f(x) %\exp(\frac{1}{4}\frac{y^{2n}}{x^{2n}})erfc(\frac{1}{2}\frac{y^n}{x^n})dx \nonumber
%\end{equation}
\end{lem}

\newpage
\begin{lem}
The following identity
\begin{equation}
\mathcal{L}_{2n}\{\mathcal{L}_n\{f(x);u\};y\} \nonumber
\end{equation}
\begin{equation}
=\frac{1}{2ny^{2n}} \int\limits_0^\infty x^{n-1}f(x)dx-\frac{\sqrt{\pi}}{4ny^{3n}}\int\limits_0^\infty x^{n-1} \frac{1}{x^{3n}} f(\frac{1}{x}) \exp(\frac{1}{4x^{2n} y^{2n}}) erfc(\frac{1}{2x^n y^n})dx, \label{lemma2}
\end{equation}
\begin{equation}
n=2^k, k=0,1,2,\ldots \nonumber
\end{equation}

is true for the $\mathcal{L}_{2n}$-transform and the $\mathcal{L}_n$-transform, provided that each member of the assertion (\ref{lemma2}) exits.
\pr
By the definitions (\ref{laplacen}),(\ref{laplace2n}) of the $\mathcal{L}_n$-transform and the $\mathcal{L}_{2n}$-transform we have
\begin{equation}
\mathcal{L}_{2n}\{\mathcal{L}_n\{f(x);u\};y\}=\int\limits_0^\infty u^{2n-1} \exp(-y^{2n} u^{2n}) \Big[\int\limits_0^\infty x^{n-1} \exp(-u^n x^n) f(x) dx\Big] du. \label{lemma2.1}
\end{equation}
Changing the order of integration, which is permissible by absolute convergence of the integrals involved, it follows from (\ref{lemma2.1}),
\begin{equation}
\mathcal{L}_{2n}\{\mathcal{L}_n\{f(x);u\};y\}=\int\limits_0^\infty x^{n-1} f(x) \Big[u^{2n-1} \exp(-y^{2n} u^{2n} - u^n x^n) du\Big] dx. \label{lemma2.2}
\end{equation}
Substituting
\begin{equation}
-y^{2n}u^{2n} - u^n x^n =-y^{2n}\Big(u^n + \frac{x}{2y^{2n}}\Big)^2 + \frac{x^{2n}}{4y^{2n}}
\end{equation}
in the inner integral on the right-hand side of (\ref{lemma2.2}) and setting
\begin{equation}
y^n(u^n + \frac{x^n}{2y^{2n}})=t
\end{equation}
we get from (\ref{lemma2.2}) the following relation:
%\begin{equation}
%\!\!\!\!\!\!\!\!\!\!\!\!\!\!\!\!\!\!\!\!\!\!\!\!\!\!\!\!\!\!\!\!\! \mathcal{L}_{2n}\{\mathcal{L}_n\{f(x);u\};y\}=\int\limits_0^\infty %x^{n-1}f(x)\exp(\frac{x^{2n}}{4y^{2n}})\times   \nonumber
%\end{equation}
%\begin{equation}
%~~~~~~~~~~~~~~~~~~~~~~~~~~~~~~~~~~ \times \Big[\frac{1}{ny^{2n}}\int\limits_{\frac{x^n}{2y^n}}^\infty t \exp(-t^2)dt - %\frac{x^n}{2ny^{3n}}\int\limits_{\frac{x^n}{2y^n}}^\infty \exp(-t^2)dt\Big]dx,
%\end{equation}
\begin{equation}
\mathcal{L}_{2n}\{\mathcal{L}_n\{f(x);u\};y\} \nonumber
\end{equation}
\begin{equation}
=\int\limits_0^\infty x^{n-1}f(x)\exp(\frac{x^{2n}}{4y^{2n}})\Big[\frac{1}{2ny^{2n}} \exp(-\frac{1}{4} \frac{x^{2n}}{y^{2n}}) - \frac{\sqrt{\pi}}{4n} \frac{x^n}{y^{3n}} erfc(\frac{x^n}{2y^n}) \Big] dx. \label{2.10}
\end{equation}

Changing the variable integration on right-hand side of (\ref{2.10}) from $x$ to $t$ according to the transformation $x=\frac{1}{t}$, we have
\begin{equation}
\mathcal{L}_{2n}\{\mathcal{L}_n\{f(x);u\};y\}=\frac{1}{2ny^{2n}} \int\limits_0^\infty x^{n-1} f(x) dx - \frac{1}{2y^{3n}} \int\limits_0^\infty t^{-2n-1} f(\frac{1}{t}) \exp(\frac{1}{4t^{2n} y^{2n}}) erfc(\frac{1}{2t^n y^n}) dt. \label{lemma2.3}
\end{equation}
The assertion (\ref{lemma2}) follows from (\ref{lemma2.3}).
\end{lem}

\begin{cor}
If the hypothesis stated in the Lemma 2.2 is satisfied then,
\begin{equation}
\mathcal{L}_{2n}\{\mathcal{L}_n\{f(x);u\};y\}=\frac{1}{2ny^{2n}} \int\limits_0^\infty x^{n-1} f(x) dx - \frac{1}{2y^{3n}} \mathcal{L}_n\{\mathcal{L}_{2n}\{\frac{1}{x^{3n}} f(\frac{1}{x}) ;u\};\frac{1}{y}\}.  \label{corollary1}
\end{equation}
\pr
Using the definition of the $\mathcal{L}_n$-transform and the $\mathcal{L}_{2n}$-transform in (\ref{lemma2.3}), we obtain the identity (\ref{corollary1}).
\end{cor}

\begin{cor}
If the hypothesis stated in the Lemma 2.2 is satisfied, then the following identity holds true for the $\mathcal{L}_n$-transform and the $\mathcal{L}_m$-transform
\begin{equation}
\mathcal{L}_n\{f(x);y\}=\frac{m}{n} \mathcal{L}_m\{f(x^{\frac{m}{n}});y^{\frac{n}{m}}\} \label{corollary2}
\end{equation}
where $m=2^{k_1}$, $n=2^{k_2}$, ~~ $k_1,k_2 \in \{0,1,2,...\}, ~~ k_1 \neq k_2$.
\pr
Using the definition (\ref{laplacen}) of the $\mathcal{L}_n$-transform we have
\begin{equation}
\mathcal{L}_n\{f(x);y\}=\int\limits_0^\infty x^{n-1} \exp(-y^n x^n) f(x) dx.  \label{corollary2.1}
\end{equation}
Setting $x^n=t^m$ on the right-hand side of (\ref{corollary2.1}) we get
\begin{equation}
\mathcal{L}_n\{f(x);y\}=\frac{m}{n} \int\limits_0^\infty t^{m-1} \exp(-y^n t^m) f(t^{\frac{m}{n}}) dt.
\end{equation}
Changing the variable on right-hand side from $t$ to $x$ according to the transformation $t=x$, we obtain (\ref{corollary2}).
\end{cor}

\begin{ex}
We show
\begin{equation}
\int\limits_0^\infty \frac{1}{x^{2n+1}} \exp[-x^{-2n} (1-\frac{1}{4} y^{-2n})] erfc(\frac{1}{2x^n y^n}) dx = \frac{1}{n} \frac{y^n}{2y^n + 1}~~ , ~~~~(n=2^k, k=0,1,2,...).
\end{equation}
\dem
If we set
\begin{equation}
f(x)=\exp(-x^{2n}) \nonumber
\end{equation}
in the formula (\ref{corollary2}), we get
\begin{equation}
\mathcal{L}_{2n}\{\mathcal{L}_n\{\exp[-x^{-2n}(1-\frac{1}{4}y^{-2n})];u\};y\}  \nonumber
\end{equation}
\begin{equation}
=\frac{1}{2ny^{2n}} \int\limits_0^\infty x^{n-1} \exp(-x^{2n}) dx - \frac{1}{2y^{3n}} \mathcal{L}_n\{\mathcal{L}_{2n}\{x^{-3n} \exp(-x^{-2n});u\};\frac{1}{y}\}. \label{example1.1}
\end{equation}
Using the formula (\ref{lemma1}) of Lemma 2.1 we have
\begin{equation}
\mathcal{L}_n\{\mathcal{L}_{2n}\{x^{-3n} \exp(-x^{2n});u\};\frac{1}{y}\}=\frac{\sqrt{\pi}}{2n} \int\limits_0^\infty x^{-2n+1}  \exp[-x^{-2n}(1-\frac{1}{4}y^{-2n})] erfc(\frac{1}{2x^ny^n})dx.  \label{example1.2}
\end{equation}
Using the relation (\ref{example1.2}) on the right-hand side of (\ref{example1.1}) we obtain
\begin{equation}
\int\limits_0^\infty \frac{1}{x^{2n+1}} \exp[-x^{-2n}(1-\frac{1}{4}y^{-2n})] erfc(\frac{1}{2x^ny^n})dx \nonumber
\end{equation}
\begin{equation}
=\frac{2y^n}{\sqrt{\pi}} \int\limits_0^\infty x^{n-1}\exp(-x^{2n})dx - \frac{4n}{\sqrt{\pi}} y^{3n} \mathcal{L}_{2n}\{\mathcal{L}_n\{\exp(-x^{2n});u\};y\}. \label{example1.3}
\end{equation}
Setting $f(x)=\exp(-x^{2n})$ in the relation (\ref{laplacen.1}) and using the known formula \cite[p.177, Entry(10)]{erdelyi1}, the Laplace transforms on the right-hand side of (\ref{example1.3}) are given by
\begin{equation}
\mathcal{L}_n\{\exp(-x^{2n});u\}=\frac{1}{n} \mathcal{L}\{\exp(-x^{2});u^{n}\}=\frac{\sqrt{\pi}}{2n} \exp(\frac{u^{2n}}{4})erfc(\frac{u^n}{2}), \label{erdelyi2}
\end{equation}
\begin{equation}
\mathcal{L}\{\exp(\frac{u}{4})erfc(\frac{\sqrt{u}}{2});y^{2n}\}=\frac{y^{-n}}{y^n + \frac{1}{2}}, \label{erdelyi}
\end{equation}
respectively, we obtain
\begin{equation}
\frac{2y^n}{\sqrt{\pi}}\int\limits_0^\infty x^{n-1} \exp(-x^{2n}) dx - \frac{4n}{\sqrt{\pi}} y^{3n} \mathcal{L}_{2n}\{\mathcal{L}_n\{\exp(-x^{2n});u\};y\}  \nonumber
\end{equation}
\begin{equation}
=\frac{2y^n}{n\sqrt{\pi}}\int\limits_0^\infty \exp(-t^2) dt - 2y^{3n} \mathcal{L}_{2n}\{\exp(\frac{u^{2n}}{4})erfc(\frac{u^n}{2});y\} \nonumber
\end{equation}
\begin{equation}
=\frac{y^n}{n} - \frac{y^{3n}}{n} \mathcal{L}\{\exp(\frac{u}{4})erfc(\frac{\sqrt{u}}{2});y^{2n}\}=\frac{1}{n} \frac{y^n}{2y^n + 1}.
\end{equation}
\end{ex}

\begin{ex}
We show
\begin{equation}
\int\limits_0^\infty \sin x ~\exp(\frac{y^2}{4x^2}) erfc(\frac{y}{2x})dx = \sqrt{\pi} y \exp(y^2) erfc(y).  \label{example2}
\end{equation}
\dem
If we set
\begin{equation}
f(x)=\sin x
\end{equation}
in the assertion (\ref{lemma1}) of Lemma 2.1 for $n=1$, we get
\begin{equation}
\mathcal{L}\{\mathcal{L}_2\{\sin x;u\};y\}=\frac{\sqrt{\pi}}{2} \int\limits_0^\infty \sin x \exp(\frac{y^2}{4x^2}) erfc(\frac{y}{2x})dx \nonumber
\end{equation}
\begin{equation}
=\frac{2}{\sqrt{\pi}} \mathcal{L}\{\mathcal{L}_2\{\sin x;u\};y\}.
\end{equation}
Now with using the above relation and the known formula \cite[p.146, Entry (21)]{erdelyi1} we have
\begin{equation}
\mathcal{L}_2\{\sin x;u\}=\sum\limits_{n=0}^\infty \frac{(-1)^n}{(2n+1)!}\mathcal{L}_2\{x^{2n+1};u\}=\frac{\sqrt{\pi}}{4u^3} \exp(-\frac{1}{4u^2}).
\end{equation}
We deduce the assertion (\ref{example2}) of Example 2.2 as follows:
\begin{equation}
\int\limits_0^\infty \sin x ~ \exp(\frac{y^2}{4x^2}) erfc(\frac{y}{2x}) dx  \nonumber
\end{equation}
\begin{equation}
=\mathcal{L}\{\frac{1}{2u^3} \exp(-\frac{1}{4u^2});y\}=\mathcal{L}\{\frac{d}{du} [\exp(-\frac{1}{4u^2})];y\} \nonumber
\end{equation}
\begin{equation}
=\sqrt{\pi} y \exp(y^2) erfc(y). \nonumber
\end{equation}
\end{ex}

\begin{ex}
We show
\begin{equation}
\int\limits_0^\infty x^{n-1} \sin(x^n) \exp(\frac{y^{2n}}{4x^{2n}}) erfc(\frac{y^n}{2x^n})dx=\frac{\sqrt{\pi}}{n} y^n \exp(y^{2n}) erfc(y^n), \label{example3}
\end{equation}
\begin{equation}
n=2^k, ~~ k=0,1,2,\ldots. \nonumber
\end{equation}
\dem
If we set
\begin{equation}
f(x)=\sin(x^n)  \nonumber
\end{equation}
in Lemma 2.1 and use the formula (\ref{lemma1}), we obtain
\begin{equation}
\frac{2n}{\sqrt{\pi}} \mathcal{L}_n\{\mathcal{L}_{2n}\{\sin(x^n);u\};y\}=\int\limits_0^\infty x^{n-1} \sin(x^n) \exp(\frac{y^{2n}}{4x^{2n}}) erfc(\frac{y^n}{2x^n}) dx.
\end{equation}
With using the relations (\ref{example2}) and (\ref{erdelyi}) and the known formula \cite[p. 146, Entry(21)]{erdelyi1},
the Laplace transforms are given by
\begin{equation}
\mathcal{L}_{2n}\{\sin(x^n);u\}=\frac{1}{2n} \mathcal{L}_2\{\sin x^{1/2};u^{2n}\}=\frac{\sqrt{\pi}}{4nu^{3n}} \exp(-\frac{1}{4} u^{-2n}),
\end{equation}
\begin{equation}
\mathcal{L}_n\{\frac{1}{2u^{3n}} \exp(-\frac{1}{4} u^{2n});y\}=\frac{1}{n}\mathcal{L}\{\frac{1}{2u^3} \exp(-\frac{1}{4} u^2); y^n\}=\frac{\sqrt{\pi}}{n} y^n \exp(y^{2n}) erfc(y^n)
\end{equation}
%\begin{equation}
%\frac{1}{n}\mathcal{L}\{\frac{1}{2u^3} \exp(-\frac{1}{4} u^2); y^n\}=\frac{\sqrt{\pi}}{n} y^n \exp(y^{2n}) erfc(y^n)
%\end{equation}
respectively. Then we obtain the relation (\ref{example3}) of Example 2.3.
\end{ex}

\begin{ex}
We show
\begin{equation}
\!\!\!\!\!\!\!\!\!\!\!\!\!\!\!\!\!\!\!\!\!\!\!\!\!\!\!\!\!\!\!\!\!
\int\limits_0^\infty \frac{1}{x} \cos(x^n) \exp(\frac{y^{2n}}{4x^{2n}}) erfc(\frac{y^n}{2x^n}) dx  \nonumber
\end{equation}
\begin{equation}
~~~~~~~~~~~~~~~~~~~~~~~~~~~~~~~~~~=\frac{4\sqrt{\pi}}{n} y^{n} (2y^{2n} +3) \exp(y^{2n}) erfc(y^n) + \frac{8}{n} y^{2n} \label{example4}
\end{equation}
\begin{equation}
~~~~~~~~~~~~~~~~~~~~~~~~~~~~~~~~~~n=2^k,~~ k=0,1,2,... \nonumber
\end{equation}
\dem
If we set
\begin{equation}
f(x)=\frac{\cos(x^n)}{x^n}
\end{equation}
in the assertion (\ref{lemma1}) of the Lemma 2.1, we get
\begin{equation}
\mathcal{L}_n\{\mathcal{L}_{2n}\{\frac{\cos(x^n)}{x^n};u\};y\}= \frac{\sqrt{\pi}}{2n} \int\limits_0^\infty x^{n-1} \frac{\cos(x^n)}{x^n} \exp(\frac{y^{2n}}{4x^{2n}}) erfc(\frac{y^n}{2x^n}) dx.  \label{example4.2}
\end{equation}
Using the equation (\ref{laplace2n.1}) we conclude
\begin{equation}
\mathcal{L}_{2n}\{\frac{\cos(x^n)}{x^n};u\}=\frac{1}{n}\mathcal{L}_2\{\frac{\cos x}{x};u^n\}=\frac{\sqrt{\pi}}{2nu^n} \exp(-\frac{1}{4u^{2n}}). \label{example4.4}
\end{equation}
Similarly by using the identity (\ref{laplacen}) and the known formulas \cite[p.129, Entry (6);p.146, Entry (21)]{erdelyi1}, we obtain
\begin{equation}
\mathcal{L}_n\{\frac{1}{u^n} \exp(-\frac{1}{4u^{2n}});y\}=\frac{1}{n} \mathcal{L}\{\frac{1}{u} \exp(-\frac{1}{4u^2});y^n\} \nonumber
\end{equation}
\begin{equation}
=\frac{4\sqrt{\pi}}{n} y^{n} (2y^{2n} +3) \exp(y^{2n}) erfc(y^n) + \frac{8}{n} y^{2n}. \label{example4.3}
\end{equation}
Substituting the result (\ref{example4.3}) into (\ref{example4.2}), we deduce the assertion (\ref{example4}).
\end{ex}
In the following lemma, we show that the second iterate of the $\mathcal{L}_{2n}$-transform in (\ref{laplace2n}) is essentially the $\mathcal{P}_{2n}$-transform defined by (\ref{potential}).

\begin{lem}
The following iteration identity,
\begin{equation}
\mathcal{L}_{2n}^2\{f(x);y\}=\mathcal{L}_{2n}\{\mathcal{L}_{2n}\{f(x);y\};z\} = \frac{1}{2n} \mathcal{P}_{2n}\{f(x);z\}, \label{lemma4}
\end{equation}
holds true, provided that the integrals involved converge absolutely.
\pr
Using the definition (\ref{laplace2n}) of the $\mathcal{L}_{2n}$-transform, we have
\begin{equation}
\!\!\!\!\!\!\!\!\!\!\!\!\!\!\!\!\!\!\!\!\!\!\!\!\!\!\!\!\!\!\!\!\! \mathcal{L}_{2n}^2\{f(x);y\}=\mathcal{L}_{2n}\{\mathcal{L}_{2n}\{f(x);y\};z\} \nonumber
\end{equation}
\begin{equation}
~~~~~~~~~~~~~~~~~~~~~~~~~~~~~~~~~~=\int\limits_0^\infty y^{2n-1} \exp(-z^{2n} y^{2n})\int\limits_0^\infty x^{2n-1} \exp(-y^{2n} x^{2n}) f(x) dx dy. \label{lemma3.1}
\end{equation}

Changing the order of integration, which is permissible by absolute convergence of the integrals involved, it follows from (\ref{lemma3.1}),
\begin{equation}
\mathcal{L}_{2n}^2 \{f(x);y\}=\mathcal{L}_{2n}\{\mathcal{L}_{2n}\{f(x);y\};z\} \nonumber
\end{equation}
\begin{equation}
=\int\limits_0^\infty x^{2n-1} f(x) \int\limits_0^\infty y^{2n-1} \exp[-y^{2n}(x^{2n} + z^{2n})]dy dx. \label{lemma3.2}
\end{equation}
Evaluating the inner integral on the right-hand side of (\ref{lemma3.2}) we have
\begin{equation}
\mathcal{L}_{2n}^2 \{f(x);y\} = \frac{1}{2n} \int\limits_0^\infty \frac{x^{2n-1}f(x)}{x^{2n} + z^{2n}} dx.
\end{equation}
With applying the definition (\ref{potential}), we deduce the identity (\ref{lemma4}) asserted by the Lemma 2.3.
\end{lem}

\begin{ex}
We illustrate the above lemma by showing
\begin{equation}
\mathcal{P}_{2n}\{\sin(z^n u^n);x\}=\frac{\pi}{n} \exp(-z^n x^n).
\end{equation}
\dem
If we set
\begin{equation}
f(x)=\sin(z^n u^n)  \nonumber
\end{equation}
in the assertion (\ref{lemma4}) of the Lemma 2.3 and then use the following identity,
\begin{equation}
\mathcal{L}_{2n}\{\sin(z^n u^n);y\}=\frac{\sqrt{\pi} z^n}{4ny^{3n}} \exp(-\frac{z^{2n}}{4y^{2n}}),  \label{example5.1}
\end{equation}
which is easily evaluable from definition of the $\mathcal{L}_{2n}$-transform (\ref{laplace2n}) and the known formula \cite[p. 146, Entry (28)]{erdelyi1}.If we apply the $\mathcal{L}_{2n}$-transform to (\ref{example5.1}) and use the relation (\ref{laplace2n.2}), we obtain
\begin{equation}
\mathcal{P}_{2n}\{\sin(z^n u^n);x\}=2n \mathcal{L}_{2n}\{\mathcal{L}_{2n}\{\sin(z^n u^n);y\};x\}=\frac{z^n \sqrt{\pi}}{2} \mathcal{L}_{2n}\{y^{-3n} \exp(-\frac{z^{2n}}{4y^{2n}});x\}.
\end{equation}
Using the identity (\ref{corollary2}) and the known formula \cite[p. 146, Entry (27)]{erdelyi1}, we find
\begin{equation}
\mathcal{P}_{2n}\{\sin(z^n u^n);x\}=\frac{z^n \sqrt{\pi}}{2} \frac{1}{2n} \mathcal{L}\{y^{-\frac{3}{2}} \exp(-\frac{z^{2n}}{4y});x^{2n}\} \nonumber
\end{equation}
\begin{equation}
=\frac{\pi}{2n} \exp(-z^n x^n).
\end{equation}
\end{ex}

\begin{ex}
We show
\begin{equation}
\mathcal{P}_{2n}\{\frac{\cos(z^n u^n)}{u^n};x\} = \frac{\pi}{2nx^n} \exp(-z^n x^n). \label{example6}
\end{equation}
\dem
If we set
\begin{equation}
f(x)=\frac{\cos(z^n u^n)}{u^n}
\end{equation}
in the assertion (\ref{lemma4}) of the Lemma 2.3, we have
\begin{equation}
\mathcal{P}_{2n}\{\frac{\cos(z^n u^n)}{u^n};x\}=2n \mathcal{L}_{2n}\{\mathcal{L}_{2n}\{\frac{\cos(z^n u^n)}{u^n};y\};x\}.   \label{example6.4}
\end{equation}
Using the following identity that could be evaluated easily from the definition of the $\mathcal{L}_{2n}$-transform and the known formula \cite[p. 158, Entry (67)]{erdelyi1},
\begin{equation}
\mathcal{L}_{2n}\{\frac{\cos(z^n u^n)}{u^n};y\} = \frac{\sqrt{\pi}}{2ny^n} \exp(-\frac{z^{2n}}{4y^{2n}}) \label{example6.3}
\end{equation}
and with using the relation (\ref{laplace2n.2}), we obtain
\begin{equation}
\mathcal{P}_{2n}\{\frac{\cos(z^n u^n)}{u^n};x\} = \sqrt{\pi} \mathcal{L}_{2n}\{y^{-n} \exp(-\frac{z^{2n}}{4y^{2n}});x\} \nonumber
\end{equation}
\begin{equation}
=\frac{\sqrt{\pi}}{2n} \mathcal{L}\{y^{-1/2} \exp(-\frac{z^{2n}}{4y});x^{2n}\}.  \label{example6.1}
\end{equation}
Making use of another known formula \cite[p.146, Entry (27)]{erdelyi1}, the Laplace transform on the right-hand side of (\ref{example6.1}) is given by
\begin{equation}
\mathcal{L}\{y^{-1/2} \exp(-\frac{z^{2n}}{4y});x^{2n}\}=\frac{\sqrt{\pi}}{x^n} \exp(-z^n x^n). \label{example6.2}
\end{equation}
Substituting the result (\ref{example6.2}) into (\ref{example6.1}), we deduce the assertion (\ref{example6}).
\end{ex}

\begin{teo}
The following Parseval-Goldstein type identities;
\begin{equation}
\int\limits_0^\infty y^{2n-1} \mathcal{L}_{2n}\{f(x);y\} \mathcal{L}_{2n}\{g(u);y\} dy = \frac{1}{2n} \int\limits_0^\infty x^{2n-1} f(x) \mathcal{P}_{2n}\{g(u);x\} dx, \label{theorem1}
\end{equation}
\begin{equation}
\int\limits_0^\infty y^{2n-1} \mathcal{L}_{2n}\{f(x);y\} \mathcal{L}_{2n}\{g(u);y\} dy = \frac{1}{2n} \int\limits_0^\infty u^{2n-1} g(u) \mathcal{P}_{2n}\{f(x);u\} du, \label{theorem2}
\end{equation}
and
\begin{equation}
\int\limits_0^\infty x^{2n-1} f(x) \mathcal{P}_{2n}\{g(u);x\} dx = \int\limits_0^\infty u^{2n-1} g(u) \mathcal{P}_{2n}\{f(x);u\} du, \label{theorem3}
\end{equation}
hold true, provided that the integrals involved converge absolutely.
\pr
We only give here the proof of Parseval-Goldstein identity (\ref{theorem1}), as the proof of the relationship (\ref{theorem2}) is similar. The relationship (\ref{theorem3}) could be obtined easily from the assertion (\ref{theorem1}) and (\ref{theorem2}). \\ By the definition (\ref{laplace2n}) of the $\mathcal{L}_{2n}$-transform, we have
\begin{equation}
\int\limits_0^\infty y^{2n-1} \mathcal{L}_{2n}\{f(x);y\} \mathcal{L}_{2n}\{g(u);y\} dy \nonumber
\end{equation}
\begin{equation}
=\int\limits_0^\infty y^{2n-1} \mathcal{L}_{2n}\{g(u);y\} \int\limits_0^\infty x^{2n-1} \exp(-x^{2n} y^{2n}) f(x) dx dy.
\end{equation}
Changing the order of integration, which is permissible by absolute convergence of the integrals involved and using the definition (\ref{laplace2n}) once again, we get
\begin{equation}
\int\limits_0^\infty y^{2n-1} \mathcal{L}_{2n}\{f(x);y\} \mathcal{L}_{2n}\{g(u);y\} dy \nonumber
\end{equation}
\begin{equation}
= \int\limits_0^\infty x^{2n-1} f(x)\int\limits_0^\infty y^{2n-1} \exp(-x^{2n} y^{2n})\mathcal{L}_{2n}\{g(u);y\} dy dx  \nonumber
\end{equation}
\begin{equation}
=\int\limits_0^\infty x^{2n-1} f(x) \mathcal{L}_{2n}\{\mathcal{L}_{2n}\{g(u);y\};x\} dx.  \label{theorem1.1}
\end{equation}
The assertion (\ref{theorem1}) could be obtained from (\ref{theorem1.1}) and the assertion (\ref{lemma4}) of Lemma 2.3.
\end{teo}

\begin{cor}
The following identity,
\begin{equation}
\mathcal{L}_{2n}\{y^{-n} \mathcal{L}_{2n}\{f(x);\frac{1}{2^{1/n} y}\};z\}=\frac{\sqrt{\pi}}{2nz^n} \mathcal{L}_n\{x^n f(x);z\}, \label{corollary3}
\end{equation}
holds true, provided that the integrals involved converge absolutely.
\pr
By setting
\begin{equation}
g(u) = \sin(z^n u^n)
\end{equation}
in the assertion (\ref{theorem1}), we obtain
\begin{equation}
\int\limits_0^\infty y^{2n-1} \mathcal{L}_{2n}\{f(x);y\} \mathcal{L}_{2n}\{\sin(z^n u^n);y\} dy  \nonumber
\end{equation}
\begin{equation}
=\frac{1}{2n} \int\limits_0^\infty x^{2n-1} f(x) \mathcal{P}_{2n}\{\sin(z^n u^n);x\} dx.
\end{equation}
Changing the variable on the left-hand side from $y$ to $t$ according to the transformation $y=\frac{1}{2^{1/n} t}$ and using the definitions (\ref{laplacen})-(\ref{laplace2n}), we find
\begin{equation}
\mathcal{L}_{2n}\{\frac{1}{t^n} \mathcal{L}_{2n}\{f(x);\frac{1}{2^{1/n} t}\};z\}= \frac{2nz^n}{\sqrt{\pi}} \int\limits_0^\infty t^{n-1} \exp(-z^{2n} t^{2n}) \mathcal{L}_{2n}\{f(x);\frac{1}{2^{1/n} t}\}dt \nonumber
\end{equation}
\begin{equation}
=\int\limits_0^\infty x^{2n-1} \exp(-z^n x^n) f(x) dx.  \label{corollary3.1}
\end{equation}
Setting $t=y$ the assertion (\ref{corollary3}) follows from (\ref{corollary3.1}).
\end{cor}

\begin{cor}
The following identity,
\begin{equation}
\int\limits_0^\infty x^{2n-1} \sin(z^n x^n) \mathcal{P}_{2n}\{g(u);x\} dx = \frac{\pi}{2n} \mathcal{L}_n\{x^n f(x);z\}, \label{corollary4}
\end{equation}
holds true for the $\mathcal{P}_{2n}$-transform and the $\mathcal{L}_n$-transform, provided that the integrals involved converge absolutely.
\pr
By setting
\begin{equation}
f(x)=\sin(z^n x^n)
\end{equation}
in the assertion (\ref{theorem1}), we obtain
\begin{equation}
\int\limits_0^\infty y^{2n-1} \mathcal{L}_{2n}\{\sin(z^n x^n);y\} \mathcal{L}_{2n}\{g(u);y\} dy = \frac{1}{2n} \int\limits_0^\infty x^{2n-1} \sin(z^n x^n) \mathcal{P}_{2n}\{g(u);x\}dx.
\end{equation}
Using the relation (\ref{example5.1}), we have
\begin{equation}
\frac{\sqrt{\pi}z^{n}}{2}\int\limits_0^\infty \frac{1}{y^{n+1}} \exp(-\frac{z^{2n}}{4y^{2n}}) \mathcal{L}_{2n}\{g(u);y\} dy= \int\limits_0^\infty x^{2n-1} \sin(z^n x^n) \mathcal{P}_{2n} \{g(u);x\}dx. \label{corollary4.1}
\end{equation}
Changing the variable on the left-hand side of (\ref{corollary4.1}) from $y$ to $t$ according to the transformation $y=\frac{1}{2^{1/n} t}$, we obtain
\begin{equation}
\int\limits_0^\infty t^{n-1} \exp(-z^{2n} t^{2n}) \mathcal{L}_{2n}\{g(u);\frac{1}{2^{1/n}t}\}dt= \frac{1}{\sqrt{\pi}z^{n}}\int\limits_0^\infty x^{2n-1} \sin(z^n x^n) \mathcal{P}_{2n} \{g(u);x\}dx.
\end{equation}
Using the definition (\ref{laplace2n}) of the $\mathcal{L}_{2n}$-transform and the definition (\ref{laplacen}) of the $\mathcal{L}_n$-transform, we have
\begin{equation}
\mathcal{L}_{2n} \{\frac{1}{t^n} \mathcal{L}_{2n}\{g(u);\frac{1}{2^{1/n}t}\};z\}= \frac{1}{\sqrt{\pi}z^{n}}\int\limits_0^\infty x^{2n-1} \sin(z^n x^n) \mathcal{P}_{2n} \{g(u);x\}dx. \nonumber
\end{equation}
Using the assertion (\ref{corollary3}) of Corollary 2.3,
\begin{equation}
\mathcal{L}_{2n} \{\frac{1}{t^n} \mathcal{L}_{2n}\{g(u);\frac{1}{2^{1/n}t}\};z\}=\frac{\sqrt{\pi}}{2nz^n}\mathcal{L}_n\{x^nf(x);z\} \nonumber
\end{equation}
\begin{equation}
=\frac{1}{\sqrt{\pi}z^{n}}\int\limits_0^\infty x^{2n-1} \sin(z^n x^n) \mathcal{P}_{2n} \{g(u);x\}dx,
\end{equation}
we obtain the assertion (\ref{corollary4}) of Corollary 2.4.
\end{cor}

\begin{cor}
The following identity,
\begin{equation}
\mathcal{L}_{2n}\{\frac{1}{y^{3n}}\mathcal{L}_{2n}\{f(x);\frac{1}{2^{1/n}y}\};z\}=\frac{\sqrt{\pi}}{n}\mathcal{L}_n\{f(x);z\}, \label{corollary5}
\end{equation}
holds true for the $\mathcal{L}_{2n}$-transform and the $\mathcal{L}_n$-transform, provided that the integrals involved converge absolutely.
\pr
By setting
\begin{equation}
g(u)=\frac{\cos(z^n u^n)}{u^n}
\end{equation}
in the assertion (\ref{theorem1}) of the Theorem 2.1, we have
\begin{equation}
\int\limits_0^\infty y^{2n-1} \mathcal{L}_{2n}\{f(x);y\}\mathcal{L}_{2n}\{\frac{\cos(z^n u^n)}{u^n};y\}dy = \frac{1}{2n}\int\limits_0^\infty x^{2n-1} f(x) \mathcal{P}_{2n} \{\frac{\cos(z^n u^n)}{u^n};x\}dx.
\end{equation}
Using the identity (\ref{example4.4}) and the identity (\ref{example6.4}), we have
\begin{equation}
\frac{2n}{\sqrt{\pi}} \int_0^\infty y^{n-1} \exp(-\frac{z^{2n}}{4y^{2n}}) \mathcal{L}_{2n}\{f(x);y\}dy=\int\limits_0^\infty x^{n-1} \exp(-z^n x^n) f(x) dx. \label{corollary5.1}
\end{equation}
Changing the variable on the left-hand side of (\ref{corollary5.1}) from $y$ to $t$ according to the transformation $y=\frac{1}{2^{1/n} t}$, we find
\begin{equation}
\frac{n}{\sqrt{\pi}}\int\limits_0^\infty t^{-n-1} \exp(-z^{2n} t^{2n}) \mathcal{L}_{2n}\{f(x);\frac{1}{2^{1/n}t}\}dt = \int\limits_0^\infty x^{n-1} \exp(-z^n x^n) f(x) dx.
\end{equation}
Using the definition of $\mathcal{L}_n$-transform (\ref{laplacen}) and setting $t=y$, we obtain the assertion (\ref{corollary5}) of Corollary 2.5.
\end{cor}

\section{Illustrative examples}

\begin{ex}
We show
\begin{equation}
\mathcal{L}_n\{x^n erfc(a^n x^n);z\} \nonumber
\end{equation}
\begin{equation}
= \frac{1}{nz^n} \Big[\frac{\sqrt{\pi}}{z^n}-\frac{1}{2a^n}\Big] -\frac{\sqrt{\pi}}{2n} \exp(\frac{z^n}{4a^{2n}})\Big[\frac{1}{z^{2n}}-\frac{1}{2a^{2n}}\Big]erfc(\frac{z^n}{2a^n}) \label{examplee}
\end{equation}
where $a \neq 0, ~ Re(a)>0, ~ n=2^k, ~k \in \mathbb{N}$.
\dem
If we set
\begin{equation}
f(x)=erfc(a^n x^n)
\end{equation}
in the assertion (\ref{corollary3}) of Corollary 2.3 we obtain
\begin{equation}
\mathcal{L}_n\{x^n erfc(a^n x^n);z\}=\frac{2nz^n}{\sqrt{\pi}}\mathcal{L}_{2n}\{y^{-n}\mathcal{L}_{2n}\{erfc(a^n x^n);\frac{1}{2^{1/n}y}\};z\}.
\end{equation}
Changing the order of integration and evaluating the inner integral on the right hand side, we have
\begin{equation}
\mathcal{L}_{2n}\{erfc(a^n x^n);\frac{1}{2^{1/n}y}\}=\int\limits_0^\infty x^{2n-1} \exp(-\frac{x^{2n}}{4y^{2n}})\int\limits_{a^n x^n}^\infty \exp(-u^2)du dx \nonumber
\end{equation}
\begin{equation}
=\int\limits_0^\infty \exp(-u^2)\int\limits_0^{u^{1/n}/a} x^{2n-1} \exp(-\frac{x^{2n}}{4y^{2n}})dxdu \nonumber
\end{equation}
\begin{equation}
=-\frac{4y^{2n}}{2n} \Big[\int\limits_0^\infty \exp \Big(-u^2(1+ \frac{1}{4a^{2n} y^{2n}})\Big)du - \int\limits_0^\infty \exp(-u^2)du \Big].
\end{equation}
Making the change of variable
\begin{equation}
u\sqrt{1+\frac{1}{4a^{2n}y^{2n}}}=t
\end{equation}
and using the linearity of the $\mathcal{L}_{2n}$-transform, we get
\begin{equation}
\mathcal{L}_{2n}\{erfc(a^n x^n);\frac{1}{2^{1/n}y}\}=2z^n \mathcal{L}_{2n}\{y^n;z\}-4a^nz^n \mathcal{L}_{2n} \{\frac{y^{2n}}{\sqrt{1+4a^{2n} y^{2n}}};z\}, \label{examplee12}
\end{equation}
where
\begin{equation}
\mathcal{L}_{2n}\{y^n;z\}=\frac{1}{2} \mathcal{L}_n\{y^{n/2};z^2\}=\frac{1}{n(z^2)^{n+\frac{n}{2}}} \Gamma(\frac{1}{2}+1)=\frac{\sqrt{\pi}}{2nz{3n}} \label{examplee11}
\end{equation}
and
\begin{equation}
\mathcal{L}_{2n}\{\frac{y^{2n}}{\sqrt{1+4a^{2n}y^{2n}}};z\}=\int\limits_0^\infty y^{2n-1} \exp(-z^{2n}y^{2n})\frac{y^{2n}}{\sqrt{1+4a^{2n}y^{2n}}}dy. \label{examplee1}
\end{equation}
To calculate the right-hand side integration in (\ref{examplee1}), we make the change of variable
\begin{equation}
1+4a^{2n}y^{2n}=u,
\end{equation}
then we get
\begin{equation}
\mathcal{L}_{2n}\{\frac{y^{2n}}{\sqrt{1+4a^{2n}y^{2n}}};z\} \nonumber
\end{equation}
\begin{equation}
=\frac{1}{32na^{4n}} \exp(\frac{z^{2n}}{4a^{2n}})\Big[\int\limits_0^1 \exp(-\frac{-z^{2n}u}{4a^{2n}})u^{1/2}du - \int\limits_0^\infty \exp(\frac{-z^{2n}u}{4a^{2n}})u^{-1/2}du \Big].
\end{equation}
Changing the variable of the integration on the right-hand side from $u$ to $x$ according to the transformation
\begin{equation}
u^{1/2} = \frac{2a^n}{z^n}x
\end{equation}
and using the integration by parts we find
\begin{equation}
\mathcal{L}_{2n}\{erfc(a^nx^n);\frac{1}{2^{1/n}y}\} \nonumber
\end{equation}
\begin{equation}
=\frac{1}{8na^{2n}z^{2n}}+\frac{\sqrt{\pi}}{8na^nz^{3n}}\exp(\frac{z^{2n}}{4a^{2n}})erfc(\frac{z^n}{2a^n}) -\frac{\sqrt{\pi}}{16na^{3n}z^n}\exp(\frac{z^{2n}}{4a^{2n}})erfc(\frac{z^n}{2a^n}). \label{examplee13}
\end{equation}
The assertion (\ref{examplee}) could be obtained by inserting expressions (\ref{examplee11}) and (\ref{examplee13}) into relation (\ref{examplee12}).
\end{ex}

\begin{ex}
We show
\begin{equation}
\mathcal{L}_n\{x^{nv}J_v(2a^nx^n);z\}=\frac{a^{nv}2^{2v}}{n\sqrt{\pi}} (z^{2n}+4a^{2n})^{-v-\frac{1}{2}} \Gamma(v+\frac{1}{2}), \label{exampleee}
\end{equation}
where $Re(v)>-\frac{1}{2},~Re(a)>0,~n=2^k,~k \in \mathbb{N}$.
\dem
If we put
\begin{equation}
f(x)=x^{nv}J_v(2a^nx^n)
\end{equation}
in the assertion (\ref{corollary5}) of Corollary 2.5. Using the known formula \cite[p. 185, Entry (27)]{erdelyi1} we obtain
\begin{equation}
\mathcal{L}_n\{x^{nv}J_v(2a^nx^n);z\} = \frac{a^{nv}2^{2v+1}}{\sqrt{\pi}}\mathcal{L}_{2n}\{y^{n(2v-1)}\exp(-4a^{2n}y^{2n});z\} \nonumber
\end{equation}
\begin{equation}
=\frac{a^{nv}2^{2v+1}}{\sqrt{\pi}} \int\limits_0^\infty y^{2nv + n -1} \exp(-y^{2n}(z^{2n} + 4a^{2n}))dy.
\end{equation}
Changing the variable of the integration on the right-hand side from $y$ to $u$ according to the transformation
\begin{equation}
y^n(z^{2n} + 4a^{2n})^{1/2} = u ,
\end{equation}
we get
\begin{equation}
\mathcal{L}_n\{x^{nv}J_v(2a^nx^n);z\}=\frac{a^{nv}2^{2v+1}}{2n\sqrt{\pi}}(z^{2n} + 4a^{2n})^{-v-\frac{1}{2}}\int\limits_0^\infty u^{2v} \exp(-u^2)du. \nonumber
\end{equation}
Making the change of variable
\begin{equation}
u^2=x
\end{equation}
and using the definition of the $\Gamma$ function, we arrive at the relation (\ref{exampleee}).

\end{ex}

\newpage
%\bibliographystyle{plain}
%\bibliography{references}

\end{document}